\documentclass{elsart}
\usepackage{natbib, amssymb}
\usepackage{url}

\newenvironment{remark}{\refstepcounter{thm}
 \bigbreak\noindent{\bf Remark \arabic{thm} }}
 {\medbreak}

\newenvironment{definition}{\refstepcounter{thm}
 \bigbreak\noindent{\bf Definition \arabic{thm} }}
 {\medbreak}

\newenvironment{example}{\refstepcounter{thm}
\bigbreak\noindent{\bf Example \arabic{thm} }}
 {\medbreak}

\newcounter{itemcounter}

\def\proof{\rm \trivlist
 \item[\hskip \labelsep{\it Proof.}]}
\def\endproof{\qed \endtrivlist}

\def\bbb#1{{\mathbb{#1}}}

\let\sem=\bf
\let\phi=\varphi
\let\rho=\varrho
\let\theta=\vartheta
\let\epsilon=\varepsilon

\def\LT{\mathop{\rm LT}\nolimits}
\def\LC{\mathop{\rm LC}\nolimits}
\def\LM{\mathop{\rm LM}\nolimits}

\def\Mat{\mathop{\rm Mat}\nolimits}

\def\Mat{\mathop{\rm Mat}\nolimits}

\def\Ord{\mathop{\tt Ord}\nolimits}
\def\Top{\mathop{\rm Top}\nolimits}

\def\dsum{\mathop{\textstyle\bigoplus}\nolimits}

\def\hom{^{\rm hom}}
\def\deh{^{\rm deh}}
\def\sat{^{\rm sat}}
\def\sweet{^{\rm sw}}
\def\Tnr{\TT^n\langle e_1,\dots,e_r\rangle}
\def\Tmnr{\TT^{m+n}\langle e_1,\dots,e_r\rangle}

\def \G {{\mathcal G}}
\def \NN {{\mathbb N}}
\def \TT {{\mathbb T}}
\def \ZZ {{\mathbb Z}}
\def\ndiv{{\not{|}\;}}

\def\degrevlex{{\tt DegRevLex}}
\def\gens{{\rm Gens}}
\def\pairs{{\rm Pairs}}
\def\satrem{{\rm SatRem}}
\def\weaksatrem{{\rm WeakSatRem}}
\def\rem{{\rm Rem}}

\def\sugar{{\rm sugar}}
\def\S{{\rm S}}

\def\lcm{\mathop{\rm lcm}\nolimits}

\def\topdeg{{\rm TopDeg}}
\let\implies=\Longrightarrow
\def\notimplies{{\;\not\!\!\implies}}

\let\epsilon=\varepsilon
\let\rho=\varrho

\def\hom{^{\rm hom}}

\def\cocoa{\mbox{\rm
     C\kern-.13em o\kern-.07 em C\kern-.13em o\kern-.15em A}}

\journal{Journal of Symbolic Computation}
\begin{document}

\begin{frontmatter}
\title{Computing Inhomogeneous Gr\"obner Bases}

\author[Genova]{A. M. Bigatti},
\ead{bigatti@dima.unige.it}
\author[Pisa]{M. Caboara},
\ead{caboara@dm.unipi.it}
\author[Genova]{L. Robbiano\corauthref{lor}}
\corauth[lor]{Corresponding author.}
\ead{robbiano@dima.unige.it}

\address[Genova]{Dipartimento di Matematica, 
Universit\`a di Genova, Italy}
\address[Pisa]{Dipartimento di Matematica 
Leonida Tonelli, Universit\`a di Pisa, Italy} 

\begin{abstract}
In this paper we describe how an idea centered on the concept of
{\it self-saturation}\/ allows several improvements in the
computation of Gr\"obner bases via Buchberger's Algorithm. In a
nutshell, the idea is to extend the advantages of computing with
homogeneous polynomials or vectors to the general case.  
When the input data are not homogeneous, we use a technique 
described in Section~\ref{Self-Saturating Buchberger's Algorithm}:
the main tool is the procedure of a 
{\it self-saturating Buchberger's Algorithm}, and
the main result is described in Theorem~\ref{mainth}. 
Another strictly related topic is treated in
Section~\ref{The Sugar Strategy} where a mathematical foundation is given
to the {\it sugar trick}\/ which is nowadays widely used in 
most of the implementations of Buchberger's Algorithm.
A special emphasis is given in Section~\ref{Single Gradings} to the case 
of a single grading, and Section~\ref{Strategies and Timings} exhibits 
some timings and indicators showing the practical merits of our approach.
\end{abstract}


\end{frontmatter}





\section*{Introduction}
\label{Introduction}

Starting from the sixties, when  implementations 
of Buchberger's famous
algorithm for computing Gr\"obner bases became practically
feasible, it has been clear that mainly three of its
steps can be optimized. They are
the minimalization of the set of critical pairs (see for instance~\cite{CKR}), 
the optimization of 
the reduction procedure (see for instance~\cite{Br}), and
the sorting used to process the critical pairs during the algorithm.
The last aspect is less important when the input polynomials 
or vectors are homogeneous and the algorithm proceeds with 
an {\it increasing degree}\/ strategy. 
But what happens if the input data are not homogeneous?

A first answer to this question was given in the late eighties. It prescribed 
to homogenize the input data, run the algorithm, and then 
dehomogenize the computed Gr\"obner basis.
This strategy is indeed quite simple and in 
many cases works fine. Its big advantage is that critical pairs are
sorted by increasing degree and after a degree is completed
the algorithm never goes back to it. The disadvantage is that 
often it computes too large a set of polynomials or vectors.

Quite soon (we are speaking of the early nineties) a 
new tool entered the game,
the {\it sugar strategy} (see~\cite{GMNRT} and 
Section~\ref{The Sugar Strategy}). In a nutshell, the idea was to 
keep the data non homogeneous, but process the 
critical pairs as if they were coming from true homogeneous data.
This goal is achieved with the help of  a
manipulated degree called {\it sugar}\/
which substitutes the true
degree. Although a complete theoretical background
was not laid out, the  idea gained popularity.
Not much later, paper~\cite{U} described an implementation 
in the computer algebra Bergman which uses
a way to improve the ordering of the critical pairs. 
Also that source was lacking a solid theoretical foundation and it 
did not gain the same popularity as the sugar strategy.

Recently,  inspired by the new development of \cocoa\ which 
will lead to the long awaited \cocoa\ 5, we decided  to explore 
some features of Buchberger's Algorithm in great detail. 
The main purpose was to give a solid theoretical background to both the 
sugar strategy and the strategy of selection of critical pairs.
We believe that we achieved both goals, so let us explain how.
After recalling more or less well-known facts about the homogenization 
process in Section~\ref{Preliminaries}, we move quickly 
to the construction of what is called the Self-Saturating 
Buchberger's Algorithm. 

To do that, in Section~\ref{Self-Saturating Buchberger's Algorithm}
we prove several properties of the 
saturation (see Proposition~\ref{propsat}), define new notions 
such as $\overline{\sigma}$-$SatGBasis$ and 
$\sigma$-$DehomBasis$, and prove
Theorem~\ref{main} where all these notions are fully compared.
With the aid of this result we define and study several 
variants of Buchberger's Algorithm,
the Weak Self-Saturating Buchberger's Algorithm and the 
Self-Saturating Buchberger's Algorithm, and finally prove the 
desired main result, Theorem~\ref{sattheorem}. 
It simply says that the computation of
a Gr\"obner basis when the input data are inhomogeneous, 
can be performed by running any 
Weak Self-Saturating Buchberger's Algorithm.
The inspiration to achieve this goal came not only 
from the above mentioned paper~\cite{U}, but also 
from the paper~\cite{BLR} where similar strategies were described
for the efficient computation of toric ideals.

It is also noteworthy to mention the fact that the variants of the
Weak Self-Saturating Buchberger's Algorithm include the usual 
Buchberger's Algorithm as well as the algorithm obtained by 
homogenizing the input data, run the algorithm, and then 
dehomogenizing the computed Gr\"obner basis.

Section~\ref{The Sugar Strategy} is devoted to give 
a solid foundation to the 
sugar strategy which, as we said, is already used in several computer 
algebra systems.  To describe it in joking mode we could say that 
the idea is to make a recipe by adding some {\em sugar}\/ to
the degree of the inhomogeneous vectors and make them {\em sweeter}\/ in
Buchberger's Algorithm.
The main result is Proposition~\ref{sugarrules}
which describes the behavior of the sugar during 
the execution of every variant 
of Buchberger's Algorithm introduced in the previous section.
With this result we can combine the tools of 
Section~\ref{Self-Saturating Buchberger's Algorithm} 
with the sugar strategy.

Section~\ref{Single Gradings} treats the case of 
a single grading and 
shows how in that situation better results can be 
achieved (see Theorem~\ref{thm:sathomog}
and its corollaries).
The current implementation in {\cocoa} deals
only with the case of the single gradings, 
shortly to be extended to the general case, and
the final Section~\ref{Strategies and Timings} 
shows its excellent behavior on a selected bunch of examples.

Of course we are aware of many algorithms which 
optimize the computation of some Gr\"obner bases
simply by going around the problem. Among many others
we could recall the Gr\"obner walk algorithm, 
the FGLM algorithm. But we want to make it clear that 
our goal here is to optimize Buchberger's Algorithm, 
not to find alternative strategies to compute Gr\"obner bases.

As a side remark we observe that 
every Self-Saturating Buchberger's Algorithm is 
fully compatible with the SlimGB 
strategies developed in \cite{Br} and with the 
Hilbert driven algorithms
(see~\cite{T} and ~\cite{CDR}).
The integration and interplay of these
approaches will be the subject of future work.  
Finally, the readers should know that the basic 
terminology is taken from the two 
books~\cite{KR1}, \cite{KR2}.

\section{Preliminaries}
\label{Preliminaries}

We assume the basic terminology and facts explained in \cite{KR2},
Section 4.3 and Tutorial 49. Some of them are explicitly recalled
for the sake of completeness, hence most of the section contains
either well-known facts or easy generalizations of well-known facts.

\subsection{Homogenization in a polynomial ring}

In this subsection we generalize the natural concept of homogenization
to the multigraded case.

We let~$K$ be a field and $P=K[x_1,\dots,x_n]$ a polynomial ring.
Then we take a matrix $W\in\Mat_{m,n}(\ZZ)$ of rank~$m\ge 1$
and new indeterminates $y_1,\dots,y_m$
called {\sem homogenizing indeterminates}.
Moreover, we equip the
polynomial \index{homogenizing indeterminate}%
ring $\overline{P}=K[y_1,\dots,y_m, x_1,\dots,x_n]$ with
the grading defined by the matrix~$\overline{W} = (I_m \mid W)$,
where~$I_m$ denotes the identity matrix of size~$m$.

Given $m$-tuples of integers $v_j=(a_{1j}, \dots, a_{mj})$, 
$j=1, \dots, s$, we consider the tuple $(c_1, \dots, c_m)$ where
$c_k = \max\{a_{k1},\dots, a_{ks}\}$, $k=1, \dots, m$, and 
call it~$\Top(v_1, \dots, v_s)$.

\begin{definition}\label{defhomog}
Let $f\in P\setminus \{0\}$ and $F\in\overline{P}$.

\begin{enumerate}

\item Write $f=c_1 t_1 + \cdots + c_s t_s$ with $c_1,\dots,
c_s \in K \setminus \{0\}$ and distinct terms
${t_1,\dots,t_s\in \TT^n}$. Then the tuple
$\Top(\deg_W(t_1), \dots, \deg_W(t_m))$
is called the {\sem top degree} of~$f$ with respect to the
grading given by~$W$ and is denoted by~$\topdeg_W(f)$.

\item  For every $j=1,\dots,s$, we let
$\deg_W(t_j)=(\tau_{1j},\dots,\tau_{mj})\in \ZZ^m$ and
let $(\mu_1, \dots, \mu_m) = \topdeg_W(f)$.
The {\sem homogenization} of~$f$
with respect to the grading given by~$W$ is the polynomial
$$
f\hom= \sum_{j=1}^s c_j\,t_j\, y_1^{\mu_1 -\tau_{1j}} \cdots
y_m^{\mu_m-\tau_{mj}}\;\in\, \overline{P}
$$
For the zero polynomial, we set $0\hom=0$.
\index{homogenization}\index{homogenization!of a polynomial}%

\item  The polynomial $F\deh=F(1,\dots,1,x_1,\dots,x_n)\in P$
is called the {\sem dehomogenization} of~$F$ with respect
to~$y_1,\dots,y_m$.
\index{dehomogenization}\index{dehomogenization!of a polynomial}%
\end{enumerate}
\end{definition}

Given an ordering $\tau$ on $\TT^n$, the monoid of power-products
in~$P$,
we want to extend it to~$\TT^{m+n}$, the monoid of power-products of
the {\it homogenization ring}~$\overline{P}$.

\begin{definition}\label{extension}
We consider a monoid ordering $\tau$ on $\TT^n$, and
the relation~$\overline{\tau}^W$   on~$\TT^{m+n}$ which is 
defined by the following rule.
Given two terms~$t_1, t_2 \in \TT^{m+n}$, we say that
${t_1 >_{\overline{\tau}^W} t_2}$  if either
$$
{\deg_{\overline{W}}(t_1) > \deg_{\overline{W}}(t_2)}
$$
or
$$
\deg_{\overline{W}}(t_1) = \deg_{\overline{W}}(t_2) \hbox{\quad and
\quad }
t_1\deh >_{\tau} t_2\deh
$$
We call $\overline{\tau}^W$ the
{\sem extension} of $\tau$   by $W$. If it is clear which grading
we are considering, we shall simply denote it by $\overline{\tau}$.
\end{definition}

We recall that the grading represented by the matrix $W$ is said to be
{\bf positive} if each column of $W$ has some non-zero entry and the
first non-zero entry  is positive.

\begin{prop}\label{exttauprop}
Let $\tau$ be a monoid ordering on $\TT^n$
and $\overline{\tau}$ its extension by~$W\!$.

\begin{enumerate}
\item The relation $\overline{\tau}$ is a
$\deg_{\overline{W}}$-compatible monoid ordering
on~$\TT^{m+n}$.

\item If $W$ is positive, the relation $\overline{\tau}$ is a
term ordering on~$\TT^{m+n}$.

\item
Let $F \in \overline{P}$ be a non-zero homogeneous polynomial.
Then there exist $s_1,\dots,s_m \in \NN$ such that
$\LT_{\overline{\sigma}}(F) = y_1^{s_1}\cdots y_m^{s_m}
  \cdot \LT_\sigma(F\deh)$.

\item If~$\tau$ is of the form $\tau=\Ord(V)$
for a non-singular matrix $V\in\Mat_n(\ZZ)$, then we have
$\overline{\tau}=
\Ord {I_m \, W^{\mathstrut} \choose 0\;\;\; V}$.

\end{enumerate}
\end{prop}

\proof
For the easy proof see~\cite{KR2}, Proposition 4.3.14 and Lemma 4.3.16.
\endproof

\begin{remark}\label{singleW}
If $\tau$ is $\deg_W$-compatible then 
$\tau =\Ord{W \choose V'}$.
Therefore we have
$$\tiny 
\overline{\tau}=
\Ord{\pmatrix{
\,\,I_m & W \cr
0 & W \cr
0 & V'
}}=
\Ord{\pmatrix{
\,\,\,I_m & W \cr
-I_m & 0 \cr 
0 & V'
}}
$$
If $m=1$ it follows that $\overline{\tau}$
is of $y_1$-DegRev type (see~\cite{KR2}, Section 4.4)
with respect to $\deg_{\overline{W}}$.
In particular, if $m=1$ and $\tau = \degrevlex$
where {\tt Deg} denotes the standard grading on $P$
it is more common to write~$\overline{P}=K[x_1,\dots,x_n, y]$ with the
homogenizing indeterminate {\em at the end}, then
we have
$\overline{\tau} = \degrevlex$
where {\tt Deg} denotes the standard grading on $\overline{P}$.
\end{remark}

\subsection{Homogenization in a free $P$-module}
\label{Homogenization in a free P-module}

In this subsection we generalize the multihomogenization procedure to
the case of free modules.

Let $r$ be a positive integer, let $F$ denote the free $P$-module
$P^r$ and let $e_1, \dots, e_r$ be the vectors of the canonical basis
of $F$.
Then let $\delta_1, \dots, \delta_r\in \ZZ^m$
and let $\overline{F}$
be the graded free $\overline{P}$-module $\overline{F} =
\dsum_{i=1}^r\overline{P}(-\delta_i)$ where the degrees of $e_1,
\dots, e_r$ are~$\delta_1, \dots, \delta_r$ respectively.
We denote by $\Tnr$ the monomodule made by the terms $t\cdot e_i\in F$
with $t\in\TT^n$,
and by~$\Tmnr$  the monomodule made by the terms $t\cdot e_i\in
\overline{F}$ with $t\in\TT^{m+n}$.
Henceforth, when we consider module orderings on $\Tnr$
we always mean module orderings which are
compatible with a monoid ordering on~$\TT^n$ (see~\cite{KR1}
Definition 1.4.17).
The following definition  is a natural generalization of
Definition~\ref{extension}.

\begin{definition}\label{sigmabar}
We consider a module ordering $\sigma$  on $\Tnr$ and
the relation~$\overline{\sigma}^W$   on~$\Tmnr$  which is
defined by the following rule.
Given $t_1e_{i}, t_2e_{j}\in \Tmnr$, we say that
${t_1e_{i} >_{\overline{\sigma}^W} t_2e_{j}}$  if either
$$
{\deg_{\overline{W}}(t_1e_{i}) > \deg_{\overline{W}}(t_2e_{j})}
$$
or
$$
\deg_{\overline{W}}(t_1e_{i}) = \deg_{\overline{W}}(t_2e_{j})
\hbox{\quad and \quad }
t_1\deh e_{i} >_{\sigma} t_2\deh e_{j}
$$
We call $\overline{\sigma}^W$ the
{\sem extension} of $\sigma$   by $W$. If it is clear which grading
we are considering, we shall simply denote it by $\overline{\sigma}$.
\end{definition}

\begin{prop}\label{prop:LTdeh}
\label{extsigmaprop}
Let $\sigma$ be a module ordering on $\Tnr$, and
let~$\overline{\sigma}$ be its extension by~$W\!$.
\begin{enumerate}
\item The relation $\!\overline{\sigma}\!$ is a
    $\deg_{\overline{W}}$-compatible module ordering
on~$\Tmnr$.

\item If $W$ is positive, then $\overline{\sigma}$ is a
module term ordering on~$\Tmnr$.

\item
Let $U\in \overline{F}$ be a homogeneous non-zero vector.
Then there exist non negative integers $s_1,\dots,s_m$ such that
$\LT_{\overline{\sigma}}(U) = y_1^{s_1}\cdots y_m^{s_m} \cdot
\LT_\sigma(U\deh)$.

\end{enumerate}
\end{prop}

\proof
It is an easy generalization of Proposition~\ref{exttauprop}.
\endproof

Analogously to Definition~\ref{defhomog} one defines the
homogenization and dehomogenization of vectors, and
with the following proposition we recall some easy results about
homogenization and dehomogenization we will need to prove
Theorem~\ref{main}.

\goodbreak

\begin{prop}\label{hom-deh}
Let $M$ be a submodule of $F$
which is generated by
vectors $v_1,\dots, v_s$, and let~$N$ be a graded submodule
of~$\overline{F}$ which is generated by homogeneous vectors
$V_1,\dots, V_t$.

\begin{enumerate}
\item We have $(M\hom)\deh = M$.
\item
The homogenization of~$M$ can be computed via the formula
$$
M\hom = \langle v\hom_1,\dots, v\hom_s\rangle :_{{}_{\overline{F}}}
(y_1\cdots y_m)^\infty
$$
\item  The dehomogenization of~$N$ can be computed via the formula
$$
N\deh = (V_1\deh, \dots,  V_t\deh)
$$
\end{enumerate}
\end{prop}

\proof
It is an obvious generalization of~\cite{KR2},
Corollaries~4.3.5.a and~4.3.8.
\endproof


\section{Self-Saturating Buchberger's Algorithm}
\label{Self-Saturating Buchberger's Algorithm}

This section starts with some properties of the saturation
and continues with the proof of the main facts (see Theorem~\ref{main})
which will eventually lead to the algorithm for computing
inhomogeneous Gr\"obner bases.
After recalling the definition of a remainder, we write the
body of Buchberger's Algorithm to help the reader spotting the
differences
when we describe some of its variants (see Theorem~\ref{sattheorem}).
The section ends with the main  Theorem~\ref{mainth}.
We keep the notation introduced before, in particular, we
let~$\sigma$ be a module ordering on $\Tnr$, and
let $\overline{\sigma}$ be its extension by~$W\!$.

\subsection{Saturation}
\label{Saturation}

\begin{definition}\label{Usat}
Let $U \in \overline{F}$ be a homogeneous vector. We denote 
$(U\deh)\hom$ by~$U\sat$ and we call it the {\sem saturation of $U$}.
Let $N$ be a graded submodule of $\overline{F}$.
We denote $ (N\deh)\hom$ by $N\sat$
and we call it the {\sem saturation of $N$}.
\end{definition}

We are going to illustrate some fundamental properties of the
saturation.
First, it is useful to recall Proposition~\ref{prop:LTdeh}.c where we
showed that there exist $s_1,\dots, s_m \in \mathbb N$ such that 
the formula $\LT_{\overline{\sigma}}(U) = 
y_1^{s_1}\cdots y_m^{s_m} \cdot \LT_\sigma(U\deh)$ 
holds true.

\begin{prop}[\bf Properties of the Saturation]
\label{propsat}
Let $v_1,\dots, v_s$ be vectors in $F$, let $U \in \overline{F}$ be a
homogeneous non-zero vector.
\begin{enumerate}
\item There exist $r_1,\dots,r_m \in \NN$ such that
$\LT_{\overline{\sigma}}(U) = y_1^{r_1}\cdots y_m^{r_m} \cdot
\LT_{\overline{\sigma}}(U\sat)$.

\item Comparing (1) with the formula
$\LT_{\overline{\sigma}}(U) = y_1^{s_1}\cdots y_m^{s_m} \cdot
\LT_\sigma(U\deh)$,
we have $r_i\le s_i$ for $i=1,\dots, m$.

\item We have $(\LT_{\overline{\sigma}}(U))\sat = \LT_\sigma(U\deh)$.

\item  If $N$ is a graded submodule
of~$\overline{F}$, then we have the following equality $N\sat = N:_{{}_{\overline{F}}} (y_1
\cdots y_m)^\infty$.

\item  If $N$ is a graded submodule
of~$\overline{F}$, then we have the following equality 
$N\sat = \langle V\sat\, |\,  V\in N,\  V \hbox{\ homogeneous}\rangle$.

\item  If $N$ is a graded submodule
of~$\overline{F}$, then we have the following equality $(N\sat )\deh = N\deh$.

\item If $M = \langle v_1,\dots, v_s\rangle$ is a submodule of $F$,
then we have  the following equality 
$M\hom = \langle v\hom_1,\dots, v\hom_s \rangle\sat$.
\end{enumerate}
\end{prop}

\proof
Condition (1) follows from the definition.  To prove (2) we denote
by~(*) the formula $\LT_{\overline{\sigma}}(U) = y_1^{s_1}\cdots
y_m^{s_m} \cdot \LT_ \sigma(U\deh)$ in
Proposition~\ref{prop:LTdeh}.(3).  We observe that $(U\sat)\deh =
U\deh$, hence, if we apply (*) to $U \sat$ we get the equality
$\LT_{\overline{\sigma}}(U\sat) = y_1^{s'_1}\cdots y_m^{s'_m} \cdot
\LT_\sigma(U\deh)$ for suitable natural numbers~$s_1', \dots, s_m'$.
Using (1) and  (*) we get $r_i+s'_i = s_i$
for ${i = 1, \dots, m}$. Condition (3) follows from condition (2).
Next we prove (4). Let $V_1, \dots, V_t$ be
homogeneous vectors which generate $N$.  Using
Proposition~\ref{hom-deh}, we deduce the following equality $N\sat = \langle V_1\sat,
\dots, V_t\sat \rangle:_{\overline{F}} (y_1 \cdots y_m)^\infty$.
It remains to show that $N:_{\overline{F}} (y_1\cdots y_m)^\infty =
\langle V_1\sat, \dots, V_t\sat \rangle:_{\overline{F}} (y_1\cdots
y_m)^\infty$.  The inclusion~$\subseteq$ is a consequence of the
obvious relation $N \subseteq \langle V_1\sat, \dots, V_t\sat
\rangle$, while the inclusion~$\supseteq$ follows from the
observation that $V_i\sat \in N:_{\overline{F}} (y_1\cdots
y_m)^\infty$ for $i = 1, \dots, t$.  Condition (5) follows from the
definition. Clearly (6) follows from~(4) and finally, to prove (7) it
suffices to combine (4) with Proposition~\ref{hom-deh}.(2).
\endproof

\begin{definition}
Let $N$ be a graded submodule of $\overline{F}$ and
let $V_1, \dots, V_t \in N$ be non-zero homogeneous vectors.

\begin{itemize}
\item[(1)] The set $\{V_1, \dots, V_t\}$ is called a
{\em $\overline{\sigma}$-SatGBasis}\/ for $N$  if  it is a $
\overline{\sigma}$-Gr\"obner basis of a
graded submodule $\tilde{N}$
of $\overline{F}$ such that $\tilde{N}\sat = N\sat$.

\item[(2)] The set $\{V_1, \dots, V_t\}$ is called a
{\em $\sigma$-DehomBasis}\/ for $N$ if
$\{V_1\deh, \dots, V_t\deh\}$ is a $\sigma$-Gr\"obner basis of $N\deh$.
\end{itemize}
\end{definition}

\begin{prop}\label{prop:dehom}
Let $M =\langle v_1, \dots, v_s\rangle$ be a  submodule of~$F$, denote by
$N$ the module $\langle v\hom_1, \dots, v\hom_s \rangle$,
and let $\{V_1, \dots, V_t\}$
be a $\sigma$-DehomBasis for $N$.
Then $\{V_1\deh, V_2\deh, \dots, V_t\deh \}$ is a
$\sigma$-Gr\"obner basis of~$M$.
\end{prop}

\proof
The claim follows from  the chain
$$N\deh = \langle v_1\hom, v_2\hom \dots, v_s\hom\rangle\deh =
\langle v_1, v_2 \dots, v_s\rangle = M$$ where the second equality
follows from Proposition~\ref{hom-deh}.(3).
\endproof

\begin{lem}\label{dehGB}
Let $N$ be a graded submodule of $\overline{F}$ and
let $V_1, \dots, V_t \in N$ be non-zero homogeneous vectors.
Then the following conditions are equivalent.
\begin{enumerate}

\item The set $\{V_1, \dots, V_t\}$ is a  $\sigma$-DehomBasis for $N$.

\item The set $\{(\LT_{\overline{\sigma}}(V_1))\sat, \dots,
  (\LT_{\overline{\sigma}}(V_t))\sat \}$ generates
  $(\LT_{\overline{\sigma}}(N))\sat$.

\end{enumerate}
\end{lem}

\proof
By Proposition~\ref{propsat}.(3) the set $\{(\LT_{\overline{\sigma}}(V_1))\sat,
\dots,
  (\LT_{\overline{\sigma}}(V_t))\sat \}$
coincides with 
$\{\LT_{{\sigma}}(V_1\deh), \dots,\LT_{{\sigma}}(V_t\deh) \}$,
and
by Proposition~\ref{propsat}.(3),(5) we have the equality
$(\LT_{\overline{\sigma}}(N))\sat =
 \langle \LT_\sigma(V\deh)\, | \, V \in N, \ V \hbox{\ \rm
homogeneous} \rangle$.
The conclusion follows immediately.
\endproof

\begin{thm}
\label{main}
Let $N$ be a graded submodule of $\overline{F}$,
let $V_1, \dots, V_t \in N$ be non-zero homogeneous vectors, and
let us consider the following conditions.

\begin{enumerate}
\item The set $\{V_1, \dots, V_t\}$ is a $\overline{\sigma}$-Gr\"obner
  basis of $N$.

\item The set $\{V_1, \dots, V_t\}$ is a $\overline{\sigma}$-SatGBasis
  for $N$.

\item The set $\{V_1, \dots, V_t \}$ is a $\sigma$-DehomBasis for $N$.

\item The set $\{ V_1\sat, \dots, V_t\sat\}$ generates $N\sat$.
\end{enumerate}
Then we have the following chain of implications.
$$ (1) \implies (2)\implies (3) \implies (4)$$
\end{thm}

\proof
The implication $(1) \implies (2)$ is obvious.

To prove $(2) \implies (3)$ let $v  \in N\deh$.
By Proposition~\ref{propsat}.(6) and the assumption,
we have $v=V\deh$ with $V  \in \tilde{N}$.
Then there exists an index $i$ such that
$\LT_{\overline{\sigma}}(V_i)\, |\, \LT_{\overline{\sigma}}(V)$.
Consequently $\LT_\sigma(V_i\deh)\, |\, \LT_\sigma(V\deh)$,
and the proof is complete.

To prove $(3) \implies (4)$ we use the equivalent
condition of  Lemma~\ref{dehGB} and proceed by contradiction.
Let $U\in N\sat$ be a
homogeneous element  with minimal
  $(\LT_{\overline{\sigma}}(U))\sat$ among the elements
in~$N\sat$ and not in $\langle V_1\sat, \dots, V_t\sat \rangle$.
We observe that  $\LT_{\overline{\sigma}}(U)\in
\langle \LT_{\overline{\sigma}}(N\sat )\rangle\subseteq
\langle \LT_{\overline{\sigma}}(N)\rangle\sat$
and therefore, by assumption,  there exists $i$
such that $(\LT_{\overline{\sigma}}(U_i))\sat$ divides
$\LT_{\overline{\sigma}}(U)$.
We deduce that, for suitable $c\in K$ and $t \in \bbb T^n$ the vector
${V = U - c\,t\,U_i}$ has the properties: $V \in N\sat$;
$V \notin \langle V_1\sat, \dots, V_t\sat \rangle$;
$\LT_{\overline{\sigma}}(V)<_{\overline{\sigma}}
\LT_{\overline{\sigma}}(U)$.
By Definition~\ref{sigmabar}, it follows that
$(\LT_{\overline{\sigma}}(V))\sat <_\sigma (\LT_{\overline{\sigma}}(U))
\sat\!$,
a contradiction.
\endproof

\goodbreak

In the next example we show that the implications of 
Theorem~\ref{main} cannot be reversed.

\begin{example}
\label{ex:viceversa}
Let $P = \mathbb Q[x,y,z]$, $\sigma = {\tt Lex}$. We use a single
homogenizing indeterminate which we call $h$ and we write $
\overline{P} = \mathbb Q[x,y,z,h]$
according to  Remark~\ref{singleW};
then $\overline{\sigma} = {\tt DegLex}$.
Let $F_1 = xh^2 - z^3$, $F_2= x^2h - y^3$, and let~$N$
be the ideal of $\overline{P}$
generated by  $\{F_1, F_2\}$. If $F_3 = y^3h^3-z^6$
it is easy to check that~${F_3 \in N}$ and
${(\LT_{\overline{\sigma}}(N))\sat =
((\LT_{\overline{\sigma}}(F_1))\sat, (\LT_{\overline{\sigma}}(F_3))
\sat)=(x, y^3)}$.
Lemma~\ref{dehGB} implies that
$\{F_1, F_3\}$ is a $\sigma$-DehomBasis for $N$; however,
it is not  a $\overline{\sigma}$-Gr\"obner basis of any
module, therefore $(3)\notimplies (2)$.
Moreover, it is easy to see that
$F_1=F_1\sat$, $F_2 = F_2\sat$, that $(F_1, F_2) = N\sat$,
but $\{F_1\deh, F_2\deh\}$ is not a~$\sigma$-Gr\"obner basis of~$N\deh$.
Therefore $(4)\notimplies (3)$.

Now let $P = \mathbb Q[x,y,z]$, $\sigma = {\tt DegRevLex}$ and let
$\overline{P} = \mathbb Q[x,y,z,h]$. In this case we have
$\overline{\sigma} = {\tt DegRevLex}$ (see Remark~\ref{singleW}).  Let
$F_1 = x^2-yh$, $F_2 = xy-zh$, let~$N$ be the ideal
of~$\overline{P}$ generated by $\{F_1, F_2\}$, and let $F_3 =
y^2h-xzh$, so that ~$F_3\sat = y^2-xz$.  Then $\{F_1, F_2, F_3\}$ is the
reduced Gr\"obner basis of~$N$, while $\{F_1, F_2, F_3\sat\}$ is the
reduced Gr\"obner basis of a module $ \tilde{N}$ such 
that~$\tilde{N}\sat = N\sat$.  Therefore $(2)\notimplies (1)$.

\end{example}

We are going to use the above results to produce a strategy for
computing Gr\"obner bases.  First, we introduce a definition.

\goodbreak

\begin{definition}
Let $\sigma$ be a module ordering on $\Tnr$, and
let~$\overline{\sigma}$ be its extension by~$W\!$.
Let $G = \{v_1, \dots, v_s\}$ be a  set of non-zero elements
in $F$ (respectively in $\overline{F}$), and
let~$u, v$ be elements  in $F$ (respectively in $\overline{F}$). Then
$u$ is said to be {\bf a remainder} of $v$ by~$G$ if it is
the output of the division algorithm applied to $v$ and $G$.
\noindent In that case we write~$u = \rem(v, G)$.
\end{definition}

A reordering of the elements of $G$ may produce different elements
which can be called $\rem(v, G)$ (see for instance~\cite{KR1}, Theorem
1.6.4), and a variant of the division algorithm, which reduces only
the leading terms, may also produce other elements.  Therefore
$\rem(v, G)$ would really be a set.  However,
for the sake of simplicity we write $u = \rem(v, G)$ instead of $u \in
\rem(v, G)$ to mean any remainder with the property that $\LT_\sigma(u)$
not divisible by $\LT_\sigma(v_i)$ for all~$v_i\in G$.

\subsection{Self Saturation}
\label{Self Saturation}

Now we write a general version of Buchberger's Algorithm.
Instead of using the stepwise description given 
in the books~\cite{KR1} and \cite{KR2},
we prefer to concentrate on the main ingredients.
In this way it will be easier for the reader to understand
the variations presented below. We recall the notion
of~S-vector  $\S(u,v)$ of~$u,v$
(see~\cite{KR1}, Definition 2.5.1). If $\LM_\sigma(u) =c_u t_u e_i$ and
$\LM_\sigma(v) = c_vt_ve_i$, then
$\S(u,v) = \frac{\lcm(t_u,t_v)}{c_u t_u}\, u
- \frac{\lcm(t_u,t_v)}{c_v t_v}\, v$.
If $U, V$ are
homogeneous vectors, some observations on the S-vector $\S(U,V)$
are
contained in~\cite{KR2}, Remark 4.5.3.

\begin{thm}[\bf Body of Buchberger's Algorithm]\label{BA}
Let $u_1, \dots, u_s$ be non-zero vectors in~$F$ (homogeneous
non-zero vectors in~$\overline{F}$) and let $M$ be the
submodule of $F$ (graded submodule of~$\overline{F}$)  generated by
$\{u_1, \dots, u_s\}$.

\begin{enumerate}
\item {\bf (Initialization)}
$\pairs=\emptyset$, the pairs;
$\gens=(u_1,\dots, u_s)$, the generators of $M$;\\
$\G=\emptyset$, the $\sigma$-Gr\"obner basis
($\overline{\sigma}$-Gr\"obner basis) of $M$ under construction.

\item {\bf (Main loop)}
While $\gens\ne\emptyset$ and $\pairs\ne\emptyset$ do

\begin{itemize}
\item[(2a)] choose $w\in\gens$ and remove it from $\gens$,\\
or a pair $(v_i, v_j)\in \pairs$, remove it from $\pairs$, and let $w
= \S(v_i,
v_j)$;
\item[(2b)] compute a remainder $v := \rem(w, \G)$;
\item[(2c)] if $v \ne 0$ add $v$ to $\G$ and the pairs $\{(v,v_i) \mid
   v_i \in\G\}$ to $\pairs$.
\end{itemize}

\item {\bf (Output)}
Return $\G$.
\end{enumerate}

\noindent This is an algorithm which returns a $\sigma$-Gr\"obner basis
($\overline{\sigma}$-Gr\"obner basis) of~$M$, whatever
choices are made in step (2a) and whatever remainder is
computed in step (2b).
\end{thm}

\begin{definition}
Let $\G$ be a tuple of homogeneous vectors in $\overline{F}$  and
$V$ a  homogeneous vector in~$\overline{F}$.
\begin{enumerate}

\item We call {\bf weak saturating remainder}\/ of $V$ with
respect to $\G$ a vector obtained in the following way.
At each step of the division algorithm, the remainder is substituted
by an element with the same saturation.
We denote it by $\weaksatrem(V, \G)$.

\item
We call {\bf saturating remainder}\/ of $V$ with
respect to $\G$, and denote it by $\satrem(V, \G)$, a vector
$(\rem(V,\G))\sat$.

\end{enumerate}
\end{definition}

Now we describe useful variants of Buchberger's Algorithm.

\begin{definition}\label{variants}
Let $U_1, \dots, U_s$ be homogeneous vectors in~$\overline{F}$ and
let $N$ be the graded submodule  of~$\overline{F}$ generated by
$\{U_1, \dots, U_s \}$. If step (2b) in Buchberger's Algorithm is
replaced~by \\
\phantom{.} \quad (2b') compute $V:= \weaksatrem(W, \G)$;\\
\noindent the procedure is called  a {\bf Weak Self-Saturating
Buchberger's Algorithm   (WeakSelfSatBA)}.
And, in~particular, if it is replaced by the following special case of (2b')
\\
\phantom{.} \quad (2b'') compute $V:= \satrem(W, \G)$;\\
the  procedure is called the {\bf
Self-Satur\-ating Buchberger's Algorithm\\  (SelfSatBA)}.
\end{definition}

A motivation for these names comes from the following result.

\begin{thm}\label{sattheorem}
Let $U_1, \dots, U_s$ be homogeneous vectors in~$\overline{F}$ and
let $N$ be the graded submodule  of~$\overline{F}$ generated by
$\{U_1, \dots, U_s \}$.
\begin{enumerate}
\item Every WeakSelfSatBA applied to $(U_1, \dots, U_s)$
computes a $\sigma$-DehomBasis for $N$.

\item SelfSatBA applied to $(U_1, \dots, U_s)$
computes a $\overline{\sigma}$-SatGBasis for $N$.
\end{enumerate}
\end{thm}

\proof
To prove (1) note that, when we substitute a
vector with another with the same saturation, the two vectors
have the same dehomogenization.
This implies that every reduction $V:= \weaksatrem(W, \G)$ mirrors
a reduction of~$W\deh$ by $\G\deh = \{U\deh \mid U \in \G\}$
with only one possible exception:
though~$V\deh$ might still be reducible by $\G\deh$,
we may choose not to substitute~$V$ with an element with the same
saturation (which would allow the ``mirror'' reduction by $\G$), we
go to step (2c) and add $V$ to $\G$. 
In this case the ``mirror'' reduction will be later
performed as a pair.  Note that since~$\LT(V)$ is not divisible by any
leading term in $\G$ this process terminates by Dickson's
Lemma, and the output is a set of vectors $\{V_1, \dots, V_t\}$ such
that $\{V_1\deh, \dots, V_t\deh\}$ is a $\sigma$-Gr\"obner basis of
$\langle U_1\deh, \dots, U_s\deh \rangle$ which is $N\deh$ by
Proposition~\ref{hom-deh}.(3).
To prove (2) we observe that all the
replacements of \rem\ with \satrem\ are equivalent to having added
some element $V\sat$ to $\gens$ and having chosen it in step (2a) just
before choosing $V$ which would consequently reduce to 0 via~$V\sat$.
\endproof

\begin{remark}
Reconsider $F_1 = xh^2 - z^3$, $F_3 = y^3h^3-z^6$ in
$\overline{P} = \mathbb Q[x,y,z,h]$, $\overline\sigma = {\tt DegLex}$
from Example~\ref{ex:viceversa} and run WeakSelfSatBA.
\begin{itemize}
\item
In (2a) we choose $W = F_1$, in (2b') we get $V = F_1$, and in (2c) we
add it to $\G$.
\item
In (2a) we choose $W = F_3$, in (2b') we get $V = F_3$, and in (2c) we
add it to $\G$ and $(F_1,F_3)$ to \pairs.
\item
In (2a) we choose $W = S(F_1,F_3) = xz^6 - y^3z^3h$, in (2b') we have
these two reduction steps:
$W_1 = W h^2 - F_1 z^6 = -y^3z^3h^3 + z^9$;
$W_2 = W_1 + F_3 z^3 = 0$
and we are done.
\end{itemize}

\noindent The output is $\{F_1, F_3\}$ which is not a
$\overline\sigma$-Gr\"obner basis (see Example~\ref{ex:viceversa}).
\end{remark}

\goodbreak
We are ready to state the main result in this section.

\begin{thm}\label{mainth}  
Let $v_1, \dots, v_s$ be non-zero vectors in  $F$, let $M$ be the
submodule  generated by the set
$\{v_1, \dots, v_s \}$, and
let  $\{V_1, V_2, \dots, V_t\}$ be the output of any {\rm WeakSelfSatBA}
applied to the set $\{ v_1\hom, v_2\hom  \dots, v_s\hom\}$.
Then the set~$\{V_1\deh, \dots, V_t\deh\}$ is 
a~$\sigma$-Gr\"obner basis of~$M$.
\end{thm}

\proof
Let $N=\langle v_1\hom, v_2\hom \dots, v_s\hom\rangle$ and let
$\{V_1,\dots, V_t\}$ be the output of a WeakSelfSatBA algorithm
applied to $\{v_1\hom, v_2\hom \dots, v_s\hom \}$.
Theorem~\ref{sattheorem} implies that the set $\{V_1,\dots, V_t\}$
is a $\sigma$-DehomBasis for $N$, i.e.~that the 
set~$\{V_1\deh, V_2\deh, \dots, V_t\deh \}$ is
a $\sigma$-Gr\"obner basis of~$N\deh$.
The conclusion follows from Proposition~\ref{prop:dehom}.
\endproof

\section{The Sugar Strategy}
\label{The Sugar Strategy}

If~we look at the variants of Buchberger's Algorithm (see Definition~\ref{variants}),
we note that they differ from the ordinary algorithm (see Theorem~\ref{BA}) only
because they allow the replacement of a vector with another one with
the same saturation. Such replacements may create vectors with a
different degree, and hence the corresponding critical pairs and
reductions have also different degree.
We observe that a reduction can also be
viewed as a special S-vector as shown in the proof
of Proposition~\ref{sweetened}, so we can concentrate on S-vectors.
The idea is that we want to {\em keep the original degree}\/ every time we
actually perform such a replacement. Now it is time to become formal.

\begin{definition}\label{companion}
Let $V$, $V'$ be homogeneous vectors in $\overline{F}$.  Then $V'$ is said to be a
{\bf companion vector} of $V$ if there exist non-negative integers
$s_1, \dots, s_m $ such that $V' = y_1^{s_1}\cdots y_m^{s_m} V$.

\end{definition}

\begin{prop}\label{sweetened}
Let be given a variant of Buchberger's Algorithm.
For each homogeneous vector~$V$ which is used during
the execution of the algorithm, there exists a unique companion
vector~$V\sweet$ (here $\sweet$ means\/ {\bf sweetened})
which obeys the following rules.
\begin{enumerate}

\item
For every  input vector $U_1, \dots, U_s$ we
have $U_i\sweet = U_i$.

\item  For every pair of  vectors $U, V$
we have $S(U,V)\sweet = S(U\sweet, V\sweet)$.

\item During the
execution of the algorithm, when a vector $V$ is substituted
by another vector~$V'$ with the property that
${(V')\sat=V\sat }$, if we have
$V\sweet =y_1^{a_1}\cdots y_m^{a_m}V\sat$,
$V' = y_1^{b_1}\cdots y_m^{b_m}V\sat $
with suitable non-negative integers
$a_1, \dots, a_m,\, b_1, \dots, b_m$, then
we have
$(V')\sweet  = y_1^{c_1}\cdots y_m^{c_m}V\sat$ where
$(c_1,\dots, c_m) = \Top\big((a_1, \dots, a_m),(b_1, \dots, b_m)\big)$.
\end{enumerate}
\end{prop}

\proof
We need to prove that for each creation of a new
vector during the execution of the algorithm,
a unique companion vector is defined.
This statement is true for the input vectors by (1) and for the
S-vectors by (2). Replacement of a vector with another one with the
same saturation is taken care by (3). Every step of reduction is of the type
$U - c\, t' \,V$ with $c \in K$, $t' \in \bbb T^{m+n}\!$.
It can be viewed as $S(U,V)$ and has been considered in (2).
\endproof

\begin{definition}\label{defsweet}
We denote  the degree
$\deg_{\overline{W}}(V\sweet)$ by ${\rm \bf sugar}(V)$,
and  we denote the degree $\deg_{\overline{W}}(S(V_i,V_j)\sweet)
= \deg_{\overline{W}}(S(V_i\sweet, V_j\sweet))$ by ${\rm \bf sugar}(V_i, V_j)$.
We say that we use the {\bf sugar strategy} if the
choice of the pairs  in step (2a) is made starting with the
lowest sugar, not the lowest degree.

\end{definition}

Elementary properties of the sugar are contained in the following proposition
which turns out to be particularly
useful for a good implementation.

\begin{prop}\label{sugarrules}
Let be given a variant of Buchberger's Algorithm and
let $U,V \in \overline{F}$ be homogeneous non-zero vectors
which are used during the execution of the algorithm.

\begin{enumerate}
\item For every $U$ we have $(U\sweet)\sat = U\sat$ and
$\sugar(U)$ is componentwise greater than or equal to $ \deg_{\overline{W}}(U)$.

\item Suppose that  $U$ is reducible by $V$, let
$\LT_{\overline{\sigma}}(U) = t'\, \LT_{\overline{\sigma}}(V)$,  with
$t' = y_1^{a_1}\cdots y_m^{a_m}\, t$, $t \in \bbb T^n$.
Let $\LC_{\overline{\sigma}}(U)= c \,\LC_{\overline{\sigma}}(V)$
and let $A=U-c\,t'\,V$ be
the result of the reduction. Then we have the equality
$${\sugar(A)= \Top\Big(\sugar(U),\ \deg_{\overline{W}}(t)+\sugar(V)\Big)}$$
\end{enumerate}

\end{prop}

\proof
Property (1) follows as an immediate consequence
of Definition~\ref{companion}, so let us
prove property (2).  Let $U\sweet =y_1^{r_1}\cdots y_m^{r_m} U$,
$V\sweet = y_1^{s_1}\cdots y_m^{s_m} V$. Then we have
\begin{eqnarray}
\LT_{\overline{\sigma}}(U\sweet) &= &y_1^{r_1}\cdots y_m^{r_m}\LT_{\overline{\sigma}}(U)=
y_1^{r_1+a_1}\cdots y_m^{r_m+a_m}\, t\, \LT_{\overline{\sigma}}(V)  \\
\LT_{\overline{\sigma}}(V\sweet) &= &y_1^{s_1}\cdots y_m^{s_m}\LT_{\overline{\sigma}}(V)
\end{eqnarray}
Moreover,
\begin{eqnarray*}
&&\sugar(A) \\
&= &\sugar(U-ct'V) = \sugar(S(U,V)) = \deg_{\overline{W}}\big(S(U,V)\sweet\big)\\ &=&
\deg_{\overline{W}}\big(S(U\sweet, V\sweet)\big)
= \deg_{\overline{W}}\Big(\lcm\big(\LT_{\overline{\sigma}}
(U\sweet), \LT_{\overline{\sigma}}(V\sweet)\big) \Big)
\end{eqnarray*}
Using formulas (1) and (2) we get
\begin{eqnarray*}
&&\lcm\Big(\LT_{\overline{\sigma}}
(U\sweet), \LT_{\overline{\sigma}}(V\sweet)\Big) \\
&= &\lcm\Big( y_1^{r_1+a_1}\cdots y_m^{r_m+a_m}\, t\, \LT_{\overline{\sigma}}(V),
y_1^{s_1}\cdots y_m^{s_m} \LT_{\overline{\sigma}}(V)\Big)  \\
&=&\lcm\Big( y_1^{r_1+a_1}\cdots y_m^{r_m+a_m}\, t\, \LT_{\overline{\sigma}}(V),
y_1^{s_1}\cdots y_m^{s_m}\, t\,  \LT_{\overline{\sigma}}(V)\Big)
\end{eqnarray*}
Consequently
\begin{eqnarray*}
&&\sugar(A) \\
&=&\deg_{\overline{W}}\Big(\lcm\Big( y_1^{r_1+a_1}
\cdots y_m^{r_m+a_m}\, t\, \LT_{\overline{\sigma}}(V),
y_1^{s_1}\cdots y_m^{s_m}\, t\,  \LT_{\overline{\sigma}}(V)\Big) \Big)\\
&=& \Top\Big((r_1+a_1, \dots, r_m+a_m) +  \deg_{\overline{W}}(t\, V), \  (s_1,\dots, s_m) +  \deg_{\overline{W}}(t\, V)  \Big)\\
&=&\Top\Big(\sugar(U),\ \deg_{\overline{W}}(t)+\sugar(V)\Big)\\
\end{eqnarray*}
where the last equality follows from formulas (1) and (2).
\endproof

\begin{example}
Consider the polynomial ring $\overline{P} = K[y_1,y_2,x_1, x_2]$ graded by
$\overline{W} = { {1\; 0\;  \, 1\; 1} \choose {0\;   1\;  \,   0\;  1}}$.
Let $U = y_1^2y_2x_1^2-y_1^3x_2$, $V = y_2x_1-x_2$. We observe that
$U$ is homogeneous of degree $(4,1)$ and~$V$ is homogeneous of degree $(1,1)$.
With $c=1$, $t'=y_1^2x_1$, $t = x_1$ we have the reduction
$A = U - y_1^2x_1V = y_1^2 x_1x_2-y_1^3x_2$
which is homogeneous of degree $(4,1)$.
Now we consider two cases.

\begin{enumerate}

\item[Case 1]
Assume that $U\sweet = U$, $V\sweet = y_1V$ so that $\sugar(U) = (4,1)$ and
$\sugar(V) = (2,1)$.
According to Proposition~\ref{sugarrules}.(2),
we have $\sugar(A) = \Top\big((4,1), (1,0) + (2,1)  \big) = (4,1)$.

\item[Case 2]
Assume instead that $U\sweet = U$, $V\sweet = y_1y_2V$
so that we have $\sugar(U) = (4,1)$ and
$\sugar(V) = (2,2)$.
Then $\LT_{\overline{\sigma}}(U\sweet) =
y_1^2y_2x_1^2$, $\LT_{\overline{\sigma}}(V\sweet) = y_1y_2^2x_1$.
Their fundamental syzygy is
$(y_2, -y_1x_1)$ whose degree is $(4,2)$. This fact is in agreement with
Proposition~\ref{sugarrules}.(2) for which
$\sugar(A) = \Top\big((4,1), (1,0) + (2,2)  \big) = (4,2)$.
It is interesting to observe that a rule of the type $\sugar(t\,V) = \deg(t) + \sugar(V)$
would have lead to $\sugar(y_1^2x_1V) = (5,2)$,  wrongly suggesting
that the sugar of $A$ should have to be $(5,2)$.

\end{enumerate}
\end{example}

\section{Single Gradings}
\label{Single Gradings}

In this section we restrict our attention to the case of positive
$\NN$-gradings i.e.~gradings defined by a row matrix $W$ with
positive entries.  Then we have a single homogenizing indeterminate
which will be called just $y$.  A first consideration in this
direction was made in Remark~\ref{singleW}, but we can say more.

\begin{thm} \label{thm:sathomog}
Let $W\in\Mat_{1,n}(\ZZ)$ be a row matrix with
positive entries and let $P$ be graded by $W$.

\begin{enumerate}
\item If $v$ is a non-zero vector in $F$, we have 
$\LT_{\overline{\sigma}}(v\hom) = \LT_{{\sigma}}(v)$.
\item If  $N$ is a graded submodule of $\overline{F}$, 
and $\G = \{V_1, V_2, \dots , V_t \}$ is a homogeneous
$\overline{\sigma}$-Gr\"obner basis of $N$, then
the set  $\{V_1\sat, V_2\sat \dots , V_t\sat\}$ is
a $\overline{\sigma}$-Gr\"obner basis of $N\sat$.
\end{enumerate}
\end{thm}

\proof
Claim (1) is clear. To prove claim (2)
we let~$V$ be a vector in $N\sat$; we need to show that
$\LT_{\overline{\sigma}}(V_i\sat)\ | \ \LT_{\overline{\sigma}}(V)$
for some $i \in \{1,\dots, t\}$.
Proposition~\ref{propsat}.(4) implies  that $y^a \cdot V \in N$ for
some $a\in\NN$. As a consequence
$y^a \cdot \LT_{\overline{\sigma}}(V) \in \LT_{\overline{\sigma}}(N)$,
hence there exists $V_i\in\G$ such that
$\LT_{\overline{\sigma}}(V_i)\ |\ y^a \cdot \LT_{\overline{\sigma}}(V)$.
Now,
${y} \ndiv \LT_{\overline{\sigma}}(V_i\sat)$ by (a) applied to $v=V_i
\deh$,
hence $\LT_{\overline{\sigma}}(V_i\sat)\ |\ \LT_{\overline{\sigma}}(V)$
and this concludes the proof.
\endproof

\begin{cor}\label{ennesat}
Let $N$ be a graded
submodule of $\overline{F}$, and
let  $\{V_1, V_2, \dots, V_t\}$ be the output of  {\rm SelfSatBA}
applied to a set of homogeneous generators of $N$.
\begin{enumerate}
\item We have $V_i = V_i\sat$ for $i =1,\dots, t$.
\item The set $ \{V_1, V_2, \dots, V_t \}$ is
a~$\overline{\sigma}$-Gr\"obner basis of~$N\sat$.
\end{enumerate}
\end{cor}

\proof
Using Thereom~\ref{sattheorem} we deduce that
$\{V_1, V_2, \dots, V_t\}$ is a $\overline{\sigma}$-Gr\"obner basis
of a graded submodule
$\tilde{N}$ of $\overline{F}$ such that $\tilde{N}\sat= N\sat$.
On the other hand, by construction SelfSatBA produces as output
saturated
vectors.
Therefore $V_i = V_i\sat$ for $i =1,\dots, t$.
Now we can use the above theorem to deduce that $\{V_1, V_2, \dots,
V_t\}$ is
a homogeneous $\overline{\sigma}$-Gr\"obner basis of~$\tilde{N}\sat = N
\sat$,
and the proof is complete.
\endproof

\begin{cor}\label{emmehom}
Let $v_1, \dots, v_s$ be non-zero vectors in  $F$, let $M$ be the
submodule  generated by
$\{v_1, \dots, v_s \}$, and
let  $\{V_1, V_2, \dots, V_t\}$ be the output of {\rm SelfSatBA}
applied to the set $\{ v_1\hom, v_2\hom  \dots, v_s\hom\}$.
Then $ \{V_1, V_2, \dots, V_t \}$ is a
$\overline{\sigma}$-Gr\"obner basis of $M\hom$.
\end{cor}

\proof
It follows from Corollary~\ref{ennesat}.(2) and Proposition~
\ref{propsat}.(7).
\endproof

\begin{example}
The following example shows that the above theorem and its corollary
cannot be extended to $\NN^m$-gradings defined
by matrices with $m>1$. The main reason is that (1)
of the above theorem is not true anymore.
Let $P = K[x_1, x_2, x_3]$ with
$$
\sigma ={\tiny
\Ord{\pmatrix{
1&1&1\cr  0&0&1\cr  0&1&0
}}
}
\hbox{\quad and \quad }
W =
{\tiny
\pmatrix{
1&1&1\cr  
1&0&1
}
}
$$
If we let
$\overline{P} = K[y_1,y_2, x_1, x_2, x_3]$,
we have
$$
\overline{W} =
{\tiny
\pmatrix{
1&0&1&1&1\cr  0&1& 1&0&1
}
}
\hbox{\quad and \quad}
\overline{\sigma} = \Ord
{\tiny 
\pmatrix{
1&0& 1&1&1 \cr  0&1& 1&0&1\cr  0&0&1&1&1\cr  0&0&0&0&1\cr 0&0&0&1&0
}
}
$$
We consider the ideal $I$ generated by
$\{ x_1x_3 - y_1y_2x_3, \  y_2x_2^2-y_1x_1\}$ in $\overline{P}$. We
check that the element
$x_2^2x_3 -y_1^2x_3$ is not in  $I$, but the element $y_2(x_2^2x_3 -y_1^2x_3)$
which is equal to
$x_3(y_2x_2^2-y_1x_1) + y_1(x_1x_3 - y_1y_2x_3)$ is in $I$ and therefore 
the element $x_2^2x_3 -y_1^2x_3$ is in $I\sat$. Consequently, if we let
$v = x_2^2-x_1$, we see that $v\hom = y_2x_2^2-y_1x_1$ and hence
$\LT_{\overline{\sigma}}(v\hom) \ne \LT_{{\sigma}}(v)$. Moreover
$\{ x_1x_3 - y_1y_2x_3, \  y_2x_2^2-y_1x_1\}$ is the
reduced $\overline{\sigma}$-Gr\"obner
basis of $I$ and both polynomials are saturated, but it cannot be the
reduced $\overline{\sigma}$-Gr\"obner basis of~$I\sat$, since we have
just seen that~${I \ne I\sat}$.
\end{example}

\section{Strategies and Timings}
\label{Strategies and Timings}

In this paper we restrict our investigation and implementation in
{\cocoa} to the case of the single grading. The implementation 
is prototypical and  it is planned to include its final form in the 
forthcoming \cocoa\ 5.

We have already mentioned that a way to compute a Gr\"obner basis
with inhomogeneous data is to homogenize the input data, compute the 
Gr\"obner basis and then dehomogenize the result.
This strategy is achieved using a Weak Self-Saturating 
Buchberger's Algorithm where the choice is to 
never saturate and choose the pair or generator 
of lowest degree in step (2a).

For degree  compatible orderings and inhomogeneous 
input, the Self-Saturating 
Buchberger's Algorithm is nothing but  the standard 
Buchberger's Algorithm with sugar. 
In step (2a) we choose the pair or
generator of lowest sugar.
The usage of homogeneous data makes the 
computation of the sugar slightly more complicated.
The result is a small overhead.

Even if we said that we {\em always} saturate, we do not need to
saturate after every reduction step, but we saturate only at the end
when the vector (or polynomial) is no longer reducible, thus avoiding
the costly operation of saturating.

\newpage

The file containing the text of the examples discussed here can be found at\\
{\tt http://cocoa.dima.unige.it/research/papers/BigCabRob09.cocoa}

The \verb+c7+ example is the classical cyclic $7$ 
system, non homogeneous.
Examples \verb+mora9+, \verb+hairer2+, 
\verb+Butcher+ and \verb+Kin1+  are well
known in the literature. Example \verb+t51+ is 
an implicitazion problem. 
Example \verb+Lex+ is a zero dimensional ideal in a polynomial ring with
three indeterminates, whose Gr\"obner basis is
computed with respect to {\tt Lex}, while \verb+Elim+ is an
elimination problem with three polynomials in five indeterminates.

\verb+A+ and \verb+H+ stay for the sugar and 
homogeneous version of the standard 
Buchberger algorithm respectively and $S$ for the 
self-saturating version.
For every example we examine some experimental 
data about the
Buchberger's Algorithm performance, namely 
cardinality of a reduced Gr\"obner
basis (before dehomogenizing in 
the \verb+H+ and \verb+S+ cases), the number 
of S-polynomials reduced and the number
of pairs considered during a run, plus the 
time spent during the computation. 
The timings are for a special version of 
the CoCoALib-0.99
on a Intel Core2Duo system with 
2MB RAM running Linux openSUSE 10.3.
$$ 
\begin{array}{ll}
\begin{array}{ll}
\\
\\
{\rm GB Len}    \\
{\rm Poly Red}  \\
{\rm Pairs Ins} \\
{\rm Time}      \\
\end{array}
      \quad
\begin{array}{rrrrr}
{\bf c7}&&{\rm  A}  &{\rm H}   &{\rm S}\\
\hline                  	          \\
	 &&209	     &443	 &209	\cr       
	 &&2060       & 2199	 &2060    \cr       
	 &&61549      &97910	 &61549  \cr       
	 &&4.52s      &3.18s	 &4.76s   \cr      
\end{array}

&\qquad

\begin{array}{rrrrr}
{\bf hairer2}& {\rm  A}  &{\rm H}   &{\rm S} \\
\hline                                      \\
	      &72       &506	   &72     \cr
	      &560      & 3149     &560     \cr
	      &6905     &127765    &6905   \cr
	      &1.26s    &16.00s    &1.40s   \cr
\end{array}
\end{array}
$$ 
$$ 
\begin{array}{ll}
\begin{array}{ll}
\\
\\
{\rm GB Len}    \\
{\rm Poly Red}  \\
{\rm Pairs Ins} \\
{\rm Time}      \\
\end{array}
\begin{array}{rrrrr}
{\bf t51P} & {\rm A} & {\rm H} &{\rm S} \\
\hline                                   \\
 	   & 6	     & 80      & 58	 \\
 	   &76	     &239      & 191	 \\
 	   &242      &3160     & 1715	 \\
	   &3.40s    &5.23s    & 0.64s   \\
\end{array}

&\qquad

\begin{array}{rrrrr}
{\bf mora9}& {\rm A}  & {\rm H}   & {\rm S}  \\
\hline                                        \\
           &2266      & 3977	   & 2552     \\
           &22099     & 43513     &  25350   \\
           &2657075   & 7906276   & 3368371  \\
           &9.91s     & 33.05s    & 12.64s   \\
\end{array}
\end{array}
$$ 

$$ 
\begin{array}{ll}
\begin{array}{ll}
\\
\\
{\rm GB Len}    \\
{\rm Poly Red}  \\
{\rm Pairs Ins} \\
{\rm Time}      \\
\end{array}
      \quad
\begin{array}{rrrrr}
{\bf Butch} &{\rm A}    &{\rm H}    &{\rm S}  \\
\hline                                            \\
	     &23	  &188         &28        \\
	     &369	  &987         &516	  \\
	     &5635	  &17578       &9256	  \\
	     &1.85s	  &3.38s       &1.40s	  \\
\end{array}

&\qquad

\begin{array}{rrrrr}
{\bf Kin1}&{\rm A} &{\rm H} &{\rm S}  \\ 	     
\hline              	        	  \\         
	   &43       & 477     & 135	  \\         
	   &625      & 4418    & 1471	  \\         
	   &10779    & 113526  & 23124    \\         
	   &3.34s    & 10.80s  & 1.57s    \\         
	   \end{array}
\end{array}
$$ 

$$ 
\begin{array}{ll}
\begin{array}{ll}
\\
\\
{\rm GB Len}    \\
{\rm Poly Red}  \\
{\rm Pairs Ins} \\
{\rm Time}      \\
\end{array}
      \quad
\begin{array}{rrrr}
{\bf Lex} &{\rm A}    &{\rm H}    &{\rm S}  \\
\hline                                            \\
	     &4	  &122         &122        \\
	     &409	  &345         &345	  \\
	     &1465	  &7381       &7381	  \\
	     &4.40s	  &0.59s       &0.61s	  \\
\end{array}

&\qquad

\begin{array}{rrrr}
{\bf Elim}&{\rm A} &{\rm H} &{\rm S}  \\ 	     
\hline              	        	  \\         
	   &99      & 353     & 353	  \\         
	   &845      & 1488    & 1488	  \\         
	   &14330    & 62128  & 62128    \\         
	   &68.80s    & 24.45s  & 24.56s    \\         
	   \end{array}
\end{array}
$$ 
The first two Gr\"obner bases are computed with respect to the  \degrevlex\
ordering; we notice that the self-saturating algorithm behavior is the same as
the standard algorithm,  the only difference being some overhead in the
saturating case, due to more complex sugar computations, as expected.

The last four Gr\"obner bases are computed with respect to lexicographic or
elimination orderings; we see that in several cases the saturating algorithm
offers an efficient alternative to the standard/homogenizing algorithm.

We also notice that the numbers entered into the tables offer a very 
partial indication of the complexity of Gr\"obner basis computations: 
only in two out of six cases the fastest
algorithm is the one with the lowest indicators. In fact, the 
possibility of using a large pool of reducers can lead to faster reductions
and hence to an overall better performance of the algorithm itself, as it 
is clearly seen in most of the H and S versions.

\begin{flushright}
\small\it
Quantity has a
quality all its own.\\
\rm (Iosif Vissarionovich Stalin)
\end{flushright}

\section*{Acknowledgements}
\vskip -1cm
This paper was inspired by the special atmosphere 
of the Castle of Hagenberg (RISC, Austria) during the ACA 2008 session on 
{\it Gr\"obner Bases and Applications}.
In the ``realm" of Bruno Buchberger one cannot refrain from trying to 
improve Buchberger's Algorithm.




\begin{thebibliography}{99}

\bibitem{BLR} A. Bigatti, R. La Scala, L. Robbiano, 
{\it Computing Toric Ideals}, J. Symbolic
Comp. {\bf 27}, 351-- 365 (1999).

\bibitem{Br}  M. Brickenstein. 
{\it Slimgb: Gr\"obner Bases with Slim Polynomials}, in:
Rhine Workshop on Computer Algebra, Proceedings of 
RWCA'06, Basel, March 2006.55--66.

\bibitem{B1} B. Buchberger, {\it Ein Algorithmus zum Auffinden
der Basiselemente des Restklassenrings nach einem
nulldimensionalen Polynomideal},
Ph.D.\ Thesis. Universit\"at Innsbruck, 1979.

\bibitem{B2} B. Buchberger, {\it A criterion for detecting unnecessary
reductions in the construction of Groebner bases},
Proc.\ EUROSAM 79, Springer LNCS {\bf 72}, (1979) 3--21.

\bibitem{B3} B. Buchberger, {\it Groebner Bases: An Algorithmic Method
in Polynomial Ideal Theory}, in: (N.K. Bose, Ed.) {\em Multidimensional
Systems Theory.} D.\ Reidel Publ.\ Comp. Pp.  (1985) 184--232.

\bibitem{CDR} M.\ Caboara, G.\ De Dominicis and L.\ Robbiano,
{\it Multigraded Hilbert Functions and Buchberger Algorithm},   
Proceedings of the ISSAC '96, (1996)72--85.

\bibitem{CKR} M.\ Caboara, M.\ Kreuzer and L.\ Robbiano,
{\it Efficiently Computing Minimal Sets of Critical Pairs},   
J. Symb. Comput.\  Vol 38, pp 1169--1190 (2004).

\bibitem{CoCoA} The {\cocoa}Team, {\it {\cocoa}: a system for doing
Computations in Commutative Algebra},
available at \url{http://cocoa.dima.unige.it}.

\bibitem{GM} Gebauer and M\"oller, {\it On an installation 
of Buchberger's algorithm}, J.\ Symbolic Computation {\bf 6} (1987), 257--286

\bibitem{GMNRT} A. Giovini, T. Mora, G. Niesi, L. Robbiano, C. Traverso,
{\it ``One sugar cube, please"
OR selection strategies in the Buchberger algorithm},
Proceedings of ISSAC,  Bonn 1991, ACM Press, S.M. Watt Editor, 49--54

\bibitem{KR1} M.\ Kreuzer and L.\ Robbiano, {\it Computational 
Commutative Algebra 1}, Springer, Heidelberg 2000.

\bibitem{KR2} M.\ Kreuzer and L.\ Robbiano, {\it Computational 
Commutative Algebra 2}, Springer, Heidelberg 2005.

\bibitem{T} C.\ Traverso, {\it Hilbert Functions and the 
Buchberger Algorithm}, 
J.\ Symbolic Computation {\bf 22} (1996), 355--376

\bibitem{U} V.\ Ufnarovski, {\it On the Cancellation Rule in the
Homogenization}
Computer Science J. of Moldova {\bf 1}, (1993).

\end{thebibliography}
\end{document}